\documentclass[11pt,a4paper,reqno,oneside]{amsart}
\usepackage[left=30mm, right=30mm, bottom=25mm, top=25mm,includefoot]{geometry}
\usepackage[noBBpl]{mathpazo}
\usepackage{graphicx, aligned-overset, amsmath, mathtools, tikz-cd, amsthm}
\usepackage[pdftex,
breaklinks=true,
bookmarksnumbered=true,
colorlinks=true,
linkcolor=ContrastRed,
citecolor=ContrastGreen,
linktocpage=true
]{hyperref}

\definecolor{MyRed}{HTML}{f0b5b4}
\definecolor{ContrastRed}{HTML}{dd5755}
\definecolor{MyGreen}{HTML}{5b8b8b}
\definecolor{ContrastGreen}{HTML}{45a1a1}

\theoremstyle{plain}
\newtheorem{theorem}{Theorem}[section]
\newtheorem{question}[theorem]{Question}
\newtheorem{conjecture}[theorem]{Conjecture}
\newtheorem{lemma}[theorem]{Lemma}

\newtheorem{proposition}[theorem]{Proposition}
\newtheorem{corollary}[theorem]{Corollary}

\theoremstyle{definition}

\newtheorem{definition}[theorem]{Definition}

\newcommand{\B}{B}
\newcommand{\PB}{P}
\newcommand{\Br}{\textrm{Brun}}

\newcommand{\Bd}{\textrm{Bd}}
\newcommand{\Z}{\textrm{Z}}
\newcommand{\Aut}{\mathrm{Aut}}

\newcommand{\triv}{{\bf 1}}
\newcommand{\defeq}{\vcentcolon=}
\newcommand{\equivv}{\equiv}
\renewcommand{\S}{S^2}
\newcommand{\AP}{\mathrm{AP}}

\makeatletter
\newsavebox{\@brx}
\newcommand{\llangle}[1][]{\savebox{\@brx}{\(\m@th{#1\langle}\)}%
\mathopen{\copy\@brx\kern-0.5\wd\@brx\usebox{\@brx}}}
\newcommand{\rrangle}[1][]{\savebox{\@brx}{\(\m@th{#1\rangle}\)}%
\mathclose{\copy\@brx\kern-0.5\wd\@brx\usebox{\@brx}}}
\makeatother

\newcommand{\highlight}[2][yellow]{\mathchoice%
  {\colorbox{#1}{$\displaystyle#2$}}%
  {\colorbox{#1}{$\textstyle#2$}}%
  {\colorbox{#1}{$\scriptstyle#2$}}%
  {\colorbox{#1}{$\scriptscriptstyle#2$}}}%

\begin{document}
\title[\resizebox{5.3in}{!}{Action of automorphisms of pure braid groups on 
homotopy groups of two-sphere}]{Action of automorphisms of pure braid groups \\ on homotopy groups of two-sphere \ \vspace{-0.2cm}}
\author{Ilya Alekseev}
\address{Ilya~Alekseev: Saint Petersburg University, 7/9 Universitetskaya nab., St. Petersburg, 199034, Russia}
\email{ilyaalekseev@yahoo.com}

\author{Vasily Ionin}
\address{Vasily Ionin: Chebyshev Laboratory, Saint Petersburg University, 29B 14th Line, Vasilyevsky Island, St. Petersburg, 199178, Russia}
\email{ionin.code@gmail.com}

\author{Mikhail Mikhailov \ \vspace{-0.4cm}}
\address{Mikhail Mikhailov: National Research University Higher School of Economics, Saint Petersburg, Russian Federation}
\email{mikhailovmikhaild@yandex.ru}


\begin{abstract}
We examine the Moore complex of the Delta-group structure related to the pure braid groups and introduced by Berrick, Cohen, Wong, and Wu. 
We prove that the cycle and the boundary groups are invariant under all automorphisms of the pure braid groups, and thereby, we extend the results of Li and Wu on the reflection automorphism. 
We conclude that there is an induced action of all automorphisms of the pure braid groups on the homotopy groups of the two-sphere. 
Besides, we compute this action for a small number of strands.
\end{abstract}

\maketitle

\ \vspace{-1.6cm}

\tableofcontents

\section{Introduction}

In the present paper, we shed more light on the connections between braid 
groups and homotopy groups of the two-sphere \(\S\). We focus on the relation 
involving Brunnian braids.\footnote{There are other interesting connections; 
see \cite{CW08} for a general exposition.}

A braid is said to be {\it Brunnian} if it becomes trivial after one removes 
any strand of the braid.
Denote by \(\PB_n\) and \(\Br_n\) the pure braid group and the group of all 
Brunnian braids, respectively, with \(n\) strands.
In~\cite{BCWW06}, Berrick, Cohen, Wong, and Wu introduce the sequence of 
homomorphisms
\begin{align}\label{BCWWSequence}
\ldots \rightarrow \Br_{n+1} \xrightarrow{\partial_{n+1}} \Br_n 
\xrightarrow{\partial_n} \Br_{n-1} \xrightarrow{\partial_{n-1}} \ldots 
\xrightarrow{\partial_3} \Br_2 \xrightarrow{\partial_2} \Br_1
\end{align}
such that the compositions \(\partial_{n}\circ \partial_{n+1}\) are trivial, 
and prove that~\(\Z_n \slash \Bd_n \cong \pi_n(\S)\) for any~\(n \ge 2\), 
where~\(\Z_n := {\rm Ker}(\partial_n)\) and~\(\Bd_n := {\rm 
Im}(\partial_{n+1})\).
We aim to investigate the above sequence from the perspective of automorphisms 
of \(\PB_n\).
Our research originates from the study of geometric transformations of 
Brunnian braids and links (see, e.g., \cite{BM21}).

In \cite{LW09}, Li and Wu construct an automorphism~\(\theta_n\) of~\(\PB_n\), 
which maps Brunnian braids to non-Brunnian ones but preserves~\(\Z_n\).
At the same time, their results imply that the reflection 
automorphism~\(\chi_n\) of \(\PB_n\) preserves both~\(\Z_n\) and~\(\Bd_n\) and 
induces the identity automorphism of~\(\Z_n \slash \Bd_n\). Inspired by these 
facts, we establish the following result.

\begin{theorem}\label{Theorem1}
For any {\normalfont \(n \ge 2\)}, both {\normalfont \(\Z_n\)} and 
{\normalfont \(\Bd_n\)} are characteristic subgroups of {\normalfont 
\(\PB_n\)}, hence all automorphisms of {\normalfont \(\PB_n\)} induce 
automorphisms of {\normalfont \(\pi_n(\S)\)}.
\end{theorem}

Given \(n \ge 2\), we denote by \(r_n \colon {\rm Aut}(\PB_n) \longrightarrow 
{\rm Aut}(\pi_n(\S))\) the resulting homomorphism of the automorphism groups.
For \(n=2\), the groups \(\PB_2\) and \(\pi_2(\S)\) are infinite cyclic, 
and~\(r_2\) is in fact an isomorphism.
Next, for \(n=3\), the group \(\Z_3\slash\Bd_3 \cong \pi_3(\S)\) is infinite 
cyclic, generated by the image of the standard braid whose closure is the 
Borromean rings.
In this case, we compute the action of \({\rm Aut}(\PB_3)\) explicitly and 
establish the following result.

\begin{proposition}\label{Proposition1}
The homomorphism {\normalfont \(r_3\)} is non-trivial, and {\normalfont \({\rm 
Im}(r_3) = \{{\rm id}, {\rm -id}\}\)}.
\end{proposition}

Finally, since~\(\pi_4(\S) \cong \pi_5(\S) \cong \mathbb{Z}/2\mathbb{Z}\), the 
homomorphisms \(r_4\) and \(r_5\) are trivial. In general, we conjecture 
that~\({\rm Im}(r_n) = \{{\rm id}, {\rm -id}\}\), so the homomorphism \(r_n\) 
is trivial if and only if the exponent of the group \(\pi_n(\S)\) equals two 
(see Conjecture~\ref{Conjecture}).

\noindent\textbf{Acknowledgments.}
Vasily Ionin is supported by ``Native towns'', a social investment program of PJSC ``Gazprom Neft''.
The work on Theorem~\ref{Theorem1}, Lemma~\ref{ReflecionExplicitFormulas}, Lemma~\ref{ReflectionCommutesWithPartial}, and Appendix~\ref{AppendixSection}
was carried out with the financial support of the Ministry of Science and Higher Education of the Russian Federation in the framework of a scientific project under agreement \texttt{No.~075-15-2025-013}.

This research started during the Summer Research Program for Undergraduates 2021, organized by the Laboratory of Combinatorial and Geometric Structures at MIPT. A sincere ‘thank you’ to Alexey Yu. Miller for drawing the wonderful pictures for us.

\section{Preliminaries}

In this section, we recall the construction of the sequence 
(\ref{BCWWSequence}) and generators of the automorphism group \({\rm 
Aut}(\PB_n)\). We refer to \cite{GJP15} for standard notions of braid theory.

For~\(n = 2\), Theorem~\ref{Theorem1} is obvious, because~\(\Z_2 = \Br_2 = 
\PB_2\) and~\(\Bd_2 = \{\triv\}\). From now on, we fix~\(n \ge 3\).
Denote by~\(\B_n\) the braid group with \(n\) strands and by \(\Sigma_n\) the 
symmetric group of degree \(n\).
We number the strands from~\(1\) to~\(n\).\footnote{Our convention agrees 
with~\cite[p.~527]{LW09}, while in~\cite[p.~298]{BCWW06}, the numbering starts 
with~\(0\).}
Let~\(\pi_\beta \in \Sigma_n\) be the permutation that a braid~\(\beta \in 
\B_n\) induces on its strands.
The assignment~\(\beta \mapsto \pi_\beta\) determines a homomorphism 
from~\(\B_n\) to~\(\Sigma_n\), whose kernel is the pure braid group~\(\PB_n\). 
Recall that
\begin{align*}
\B_n \cong \langle \sigma_1, \ldots, \sigma_{n-1} \mid 
\sigma_k\sigma_{k+1}\sigma_k=\sigma_{k+1}\sigma_k\sigma_{k+1}; \; 
\sigma_i\sigma_j =\sigma_j \sigma_i, \; |i-j| \ge 2 \rangle
\end{align*}
and that the braids \(\{A_{i,j} \mid 1 \le i < j \le n\}\) generate~\(\PB_n\), 
where\footnote{Note that in \cite{BCWW06}, \(A_{i,j}\) is denoted by 
\(A_{i-1,j-1}\), while the convention of~\cite{LW09} agrees with ours.}
\begin{align*}
A_{i,j} &\defeq \sigma_{j-1} \sigma_{j-2} \ldots 
\sigma_{i+1}\sigma_i^2\sigma_{i+1}^{-1} \ldots 
\sigma_{j-2}^{-1}\sigma_{j-1}^{-1}.
\end{align*}
We picture these braids as shown in Figure~\ref{Figure1}, where rectangular 
boxes with small round holes stand for trivial braids.

Recall that the abelianization \(\PB_n \slash [\PB_n, \PB_n]\) is free abelian 
of rank~\(n(n-1)/2\), generated by the images of \(A_{i,j}\) (see~\cite[Theorem 
11]{GJP15}).
Roughly speaking, the abelianization homomorphism \(\PB_n \to \PB_n \slash 
[\PB_n, \PB_n]\) assigns to a pure braid, given by a word in the alphabet 
\(\{A_{i,j}^{\pm 1} \mid 1 \le i < j \le n\}\), the tuple of exponent sums of 
the letters~\(A_{i,j}\).

\subsection{Brunnian braid groups}

Given \(k \in \{1,2,\ldots,n\}\), denote by \(d_k \colon \PB_n \to \PB_{n-1}\) 
the removal of the \(k\)th strand homomorphism.\footnote{In 
\cite[p.~306]{BCWW06} and~\cite[p.~527]{LW09}, \(d_k\) is denoted 
by~\(d_{k-1}\).}
By definition,
\begin{align*}
\Br_n = \bigcap_{k=1}^n {\rm Ker}(d_k).
\end{align*}
It is easy to see that
\begin{align}\label{DkGenerators}
{\rm Ker}(d_k) = \langle A_{1,k}, A_{2,k}, \ldots, A_{k-1,k}, 
A_{k,k+1},\ldots,A_{k,n}\rangle.
\end{align}
In particular, for each braid \(\beta \in \Br_n\) and for each letter 
\(A_{i,j}\), the braid \(\beta\) can be presented as a word in the alphabet 
\(\{A_{s,t}^{\pm 1} \mid 1 \le s < t \le n\} \backslash \{A_{i,j}^{\pm 1}\}\), 
thus, the exponent sum of \(A_{i,j}\) in this record is zero. Therefore,
\begin{align}\label{BrunnIsInCommutatorSubgroup}
\Br_n \subseteq [\PB_n, \PB_n].
\end{align}

\begin{figure}
\centering
\includegraphics[height = 4.5cm]{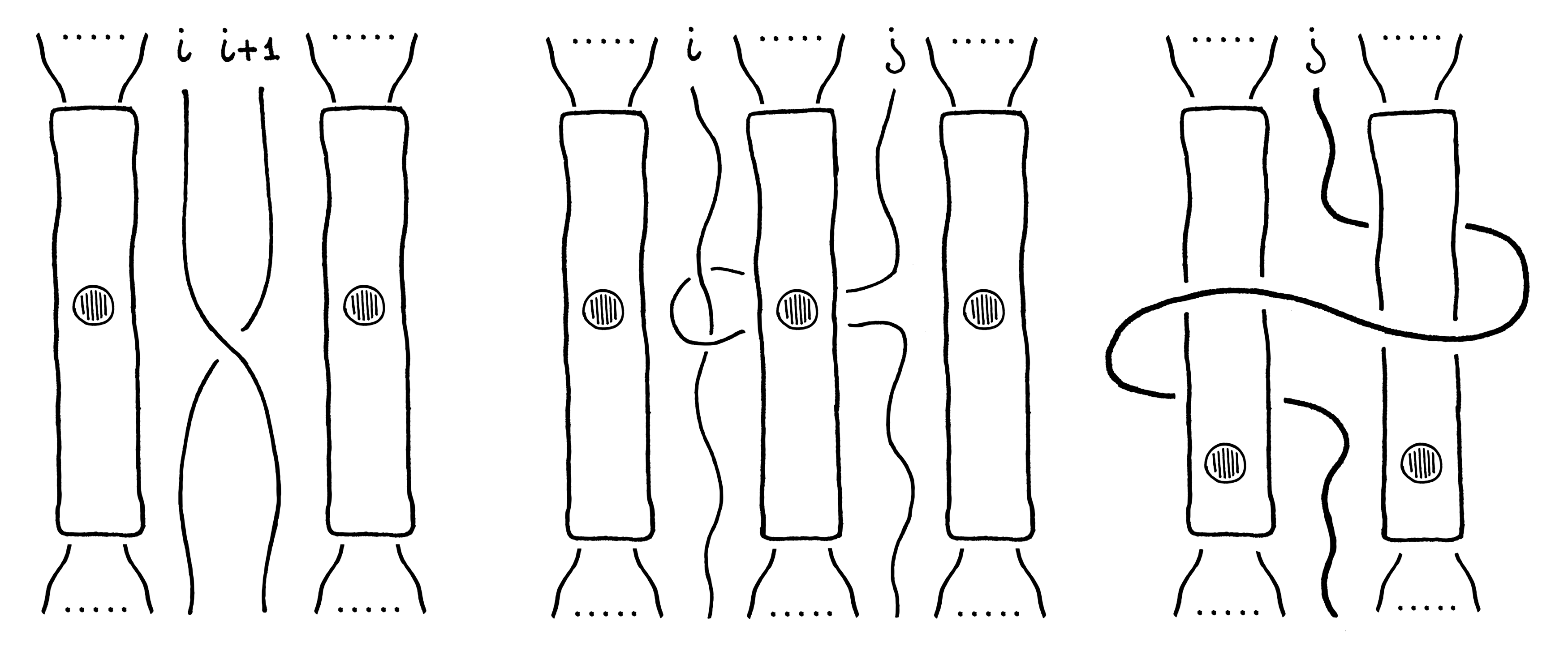}
\caption[Pictures of geometric braids \(\sigma_i\), \(A_{i,j}\), and 
\(A_{0,j}\).]{Pictures of geometric braids \(\sigma_i\), \(A_{i,j}\), and 
\(A_{0,j}\).}
\label{Figure1}
\end{figure}

\noindent Given \(j \in \{1,2,\ldots,n\}\), let\footnote{Note that in 
\cite{BCWW06}, \(A_{0,j}\) is denoted by \(A_{-1,j-1}\), while the convention 
of~\cite{LW09} agrees with ours; the fourth line in \cite[p.~530]{LW09} seems 
to be a typo.}
\begin{align*}
A_{0,j} \defeq A_{j,n}^{-1} A_{j,n-1}^{-1} \ldots A_{j,j+1}^{-1} 
A_{j-1,j}^{-1} A_{j-2,j}^{-1} \ldots A_{1,j}^{-1}.
\end{align*}
As shown in Figure~\ref{Figure1}, in the braid~\(A_{0,j}\), the~\(j\)th strand 
is wrapped around all the others.
Given indices~\(i,j \in \{0,1,\ldots,n\}\) such that~\(i < j\), let~\(A_{j,i} 
\defeq A_{i,j}\) and let~\({A_{i,i} = \triv}\) be the trivial braid. In this 
notation,
\begin{align*}%
A_{0,j} = (A_{1,j} A_{2,j} \ldots A_{n,j})^{-1}.
\end{align*}

\noindent Following~\cite{LW09}, define an automorphism
\(\theta_n \colon \PB_n \to \PB_n\) by
\begin{align*}
\theta_n(A_{i,j}) \defeq
\begin{cases}
A_{i,j}, & 2 \le i < j \le n; \\
A_{1,j}^{-1}A_{0,j} A_{1,j}, & 1 = i < j \le n.
\end{cases}
\end{align*}
See~\cite[p.~522]{LW09} for a geometric interpretation of~\(\theta_n\).

Define a homomorphism \(\partial_n \colon \PB_n \to \PB_{n-1}\) by 
\(\partial_n \defeq d_1 \circ \theta_n\). In~\eqref{BCWWSequence}, we denote 
the restriction of the homomorphism~\(\partial_n\) to~\(\Br_n\) by the same 
symbol. Let
\[\Z_n \defeq \Br_n \cap {\rm Ker}(\partial_n)\] and \[\Bd_n \defeq 
\partial_{n+1}(\Br_{n+1}).\]

The following result is a key element in the connection between brunnian 
braids and homotopy groups of the two-sphere.

\begin{lemma}\label{KeyIsomorphism}
The group {\normalfont \(\Bd_n\)} is a normal subgroup of {\normalfont 
\(\Z_n\)}, and {\normalfont\(\Z_n \slash \Bd_n \cong \pi_n(\S)\)}.
\end{lemma}

The sequence~(\ref{BCWWSequence}) is the Moore complex of a particular 
Delta-group (see Definition~\ref{Delta-groupsDef}), and its ‘homology groups’ 
are thus isomorphic to the homotopy groups of~\(\S\).
Lemma~\ref{KeyIsomorphism} is deduced in~\cite{BCWW06} and~\cite{LW09}. 
However, we found several issues in these papers.
In the Appendix, we discuss these issues and show that they do not affect the 
validity of Lemma~\ref{KeyIsomorphism} or the results of the present paper.

\subsection{Automorphisms of pure braid groups}

Recall that the center \(Z(\PB_n)\) of~\(\PB_n\) is infinite cyclic (see, 
e.g., \cite[Theorem~15]{GJP15}), generated by the full-twist braid
\begin{align*}
z_n \defeq A_{1,2} (A_{1,3} A_{2,3}) \ldots (A_{1,n} A_{2,n} \ldots A_{n-1,n}).
\end{align*}

An automorphism \(f \colon \PB_n \to \PB_n\) is said to be {\it central} if 
\(f\) induces the identity automorphism of the quotient group \(\PB_n \slash 
Z(\PB_n)\).
In this case, for any~\(x \in \PB_n\) there exists~\(m \in \mathbb{Z}\) such 
that~\(f(x) = x z_n^m\).
Denote by \({\rm Aut}_c(\PB_n)\) the subgroup of all central automorphisms 
of~\(\PB_n\).

Following \cite{BNS18}, define an automorphism \(w_n \colon \PB_n \to \PB_n\) 
by
\begin{align*}
w_n(A_{i,j}) \defeq
\begin{cases}
A_{i,j}, & 1 \le i < j < n; \\
(A_{1,n}A_{1,2}A_{1,3} \ldots A_{1,n-1})^{-1}, & (i,j) = (1,n); \\
(A_{2,n}A_{1,2}A_{2,3}\ldots A_{2,n-1})^{-1}, & (i,j) = (2,n); \\
(A_{i,n} A_{1,i} \ldots A_{i-1,i} A_{i,i+1} \ldots A_{i,n-1})^{-1}, & 3 \le i 
< j = n.
\end{cases}
\end{align*}
It is easy to check that
\begin{align*}
w_n(A_{i,j}) =
\begin{cases}
A_{i,j}, & 1 \le i < j < n; \\
A_{i,n} A_{0,i} A_{i,n}^{-1}, & 1 \le i < j = n.
\end{cases}
\end{align*}
See~\cite[p.~5]{BNS18} for a geometric interpretation of~\(w_n\) in terms of 
the mapping class group of a punctured sphere.

Let~\(\chi_n\) be the reflection automorphism of \(\PB_n\).
We denote by~\({\rm Aut}(\B_n)\) the subgroup of~\({\rm Aut}(\PB_n)\) 
generated by~\(\chi_n\) and the conjugations by all elements of~\(\B_n\).

In \cite{BNS18}, Bardakov, Neshchadim, and Singh prove that the group~\({\rm 
Aut}(\PB_n)\) is generated by the union
\begin{align}\label{AutomorphismClasses}
{\rm Aut}(\B_n) \cup {\Aut}_c(\PB_n) \cup \{w_n\}.
\end{align}
To show that both \(\Z_n\) and \(\Bd_n\) are characteristic subgroups of 
\(\PB_n\) (Theorem~\ref{Theorem1}), we prove their invariance under these 
automorphisms.

\section{Expressions for \texorpdfstring{\(\Br_n\)}{Lg} and 
\texorpdfstring{\(\Bd_n\)}{Lg}}

Given elements \(x,y \in \PB_n\), let~\([x,y] \defeq x^{-1}y^{-1}xy\).
Next, given subgroups \(A,B \subseteq \PB_n\), let~\(AB \defeq \langle xy \mid 
x \in A, \ y \in B\rangle\) and~\([A,B] \defeq \langle [x,y] \mid x \in A, \ y 
\in B\rangle\).
Given normal subgroups~\(R_1, R_2, \ldots, R_m\) of \(\PB_n\), let
\begin{align*}
[R_1, R_2, \ldots, R_m]_S \defeq \prod\limits_{\tau \in \Sigma_m} 
[[[R_{\tau(1)},R_{\tau(2)}], R_{\tau(3)}] \ldots, R_{\tau(m)}].
\end{align*}
Finally, denote by~\(\llangle \beta \rrangle\) the normal closure of a pure 
braid~\(\beta\) in~\(\PB_n\).

In \cite[Theorem~1.1]{BMVW12}, Bardakov, Mikhailov, Vershinin, and Wu 
establish the following equality for all~\(n \ge 3\):
\begin{align}\label{BrunnianSymmetricCommutator}
\Br_n = [\llangle A_{1,n} \rrangle, \llangle A_{2,n} \rrangle, \ldots, 
\llangle A_{n-1,n} \rrangle]_S.
\end{align}

Below, we deduce a similar expression for \(\Bd_n\).
For this, we compute the action of inner automorphisms of~\(\B_n\) on the set 
\begin{align*}
\{\llangle A_{i,j} \rrangle \mid 1 \le i < j \le n\}.
\end{align*}

Given \(\beta \in \B_n\), let \(\Psi_\beta \colon \PB_n \to \PB_n\) be the 
conjugation by \(\beta\), that is, \(\Psi_\beta(x) = \beta^{-1} x \beta\) for 
all~\(x \in \PB_n\).\footnote{Note that in~\cite[p.~539]{LW09}, 
\(\chi_{\sigma_k}\) denotes~\(\Psi_{\sigma_k^{-1}}\).}
It is easy to see that
\begin{align}\label{PermutationOfDs}
\Psi_\beta({\rm Ker}(d_k)) = {\rm Ker}(d_{\pi_\beta(k)}).
\end{align}

\begin{lemma}
For any \(i,j \in \{1,2,\ldots,n\}\) such that \(i < j\), we have
{\normalfont
\begin{align}\label{IntersectionOfDs}
\llangle A_{i,j} \rrangle = {\rm Ker}(d_i) \cap {\rm Ker}(d_j).
\end{align}}
\end{lemma}
\begin{proof}
By \cite[Lemma 3.8]{BMVW12}, the statement holds for \(j = n\).
For the general case, we set \(\beta = \sigma_{n-1}\sigma_{n-2}\ldots 
\sigma_j.\) Note that \(\Psi_\beta(A_{i,n}) = A_{i,j}\).
We have
\begin{align*}
{\rm Ker}(d_i) \cap {\rm Ker}(d_j) &= \Psi_\beta\left({\rm Ker}(d_i)\right) 
\cap \Psi_\beta\left({\rm Ker}(d_n)\right) = \Psi_\beta\left({\rm Ker}(d_i) 
\cap {\rm Ker}(d_n)\right) \\
&= \Psi_\beta\left(\llangle A_{i,n} \rrangle\right) = \llangle 
\Psi_\beta(A_{i,n}) \rrangle = \llangle A_{i,j} \rrangle.\qedhere
\end{align*}
\end{proof}
\medskip

\begin{corollary}
For any \(\beta \in \B_n\) and for any \(i,j \in \{1,2,\ldots,n\}\) such that 
\(i < j\), the following equality holds:
\begin{align}\label{ConjugationFormulas}
\Psi_\beta\left(\llangle A_{i,j} \rrangle\right) = \llangle 
A_{\pi_\beta(i),\pi_\beta(j)} \rrangle.
\end{align}
\end{corollary}
\begin{proof}
Indeed,
\begin{align*}
\Psi_\beta(\llangle A_{i,j} \rrangle) \overset{(\ref{IntersectionOfDs})}&{=} 
\Psi_\beta({\rm Ker}(d_i) \cap {\rm Ker}(d_j)) = \Psi_\beta({\rm Ker}(d_i)) 
\cap \Psi_\beta({\rm Ker}(d_j)) \\
\overset{(\ref{PermutationOfDs})}&{=} {\rm Ker}(d_{\pi_\beta(i)}) \cap {\rm 
Ker}(d_{\pi_\beta(j)}) \overset{(\ref{IntersectionOfDs})}{=} \llangle 
A_{\pi_\beta(i),\pi_\beta(j)} \rrangle.\qedhere
\end{align*}
\end{proof}
\medskip

Next, we provide an expression for \(\Bd_n\).

\begin{proposition}\label{BdExpressionProposition}
The following equality holds:
\begin{align}\label{BoundaryBrunnianSymmetricCommutator}
{\normalfont \Bd_n} = [\llangle A_{0,1} \rrangle, \llangle A_{0,2} \rrangle, 
\ldots, \llangle A_{0,n} \rrangle]_S.
\end{align}
\end{proposition}
\begin{proof}
Given \(m \ge 3\), let \(\beta \defeq (\sigma_{m-1}\sigma_{m-2}\ldots 
\sigma_2)\sigma_1(\sigma_{2}^{-1}\sigma_{3}^{-1}\ldots \sigma_{m-1}^{-1}) \in 
\B_m.\)
It follows from~(\ref{PermutationOfDs}) that \(\Psi_\beta\) 
preserves~\(\Br_m\). Note that the permutation \(\pi_\beta\) is a transposition 
that swaps~\(1\) and~\(m\).
We have
\begin{align*}
\Br_m = \Psi_\beta(\Br_m) \overset{(\ref{BrunnianSymmetricCommutator})}&{=} 
\Psi_\beta\left([\llangle A_{1,m} \rrangle, \llangle A_{2,m} \rrangle, \ldots, 
\llangle A_{m-1,m} \rrangle]_S\right) \\
&= [\llangle \Psi_\beta(A_{1,m}) \rrangle, \llangle \Psi_\beta(A_{2,m}) 
\rrangle, \ldots, \llangle \Psi_\beta(A_{m-1,m}) \rrangle]_S \\
\overset{(\ref{ConjugationFormulas})}&{=} [\llangle A_{1,m} \rrangle, \llangle 
A_{2,1} \rrangle, \ldots, \llangle A_{m-1,1} \rrangle]_S \\
&= [\llangle A_{2,1} \rrangle, \llangle A_{3,1} \rrangle, \ldots, \llangle 
A_{m,1} \rrangle]_S.
\end{align*}
Next, we have
\begin{align*}
\Bd_n = \partial_{n+1}(\Br_{n+1}) &= \partial_{n+1}\left([\llangle A_{2,1} 
\rrangle, \llangle A_{3,1} \rrangle, \ldots, \llangle A_{n+1,1} \rrangle]_S 
\right) \\
&= [\llangle \partial_{n+1}(A_{2,1}) \rrangle, \llangle 
\partial_{n+1}(A_{3,1}) \rrangle, \ldots, \llangle \partial_{n+1}(A_{n+1,1}) 
\rrangle]_S \\
&=[\llangle A_{0,1} \rrangle, \llangle A_{0,2} \rrangle, \ldots, \llangle 
A_{0,n} \rrangle]_S.\qedhere
\end{align*}
\end{proof}

\section{Proof of Theorem \ref{Theorem1}}

We examine the three classes of automorphisms described in 
(\ref{AutomorphismClasses}) separately.

Firstly, we prove that both \(\Z_n\) and \(\Bd_n\) are invariant under \({\rm 
Aut}_c(\PB_n)\). For this, we prove that \({\rm Aut}_c(\PB_n)\) acts 
identically on these subgroups.

\begin{lemma}\label{TrivialityOfAction}
The group {\normalfont \({\rm Aut}_c(\PB_n)\)} acts identically on 
{\normalfont \([\PB_n,\PB_n]\)}.
\end{lemma}
\begin{proof}
Let \(f\) be a central automorphism of \(\PB_n\). Then for any \(x,y \in 
\PB_n\) there exist \(s,t \in \mathbb{Z}\) such that \(f(x) = x z_n^s\) and 
\(f(y) = y z_n^t\). We have
\begin{align*}
f([x,y]) = [f(x),f(y)] = [x z_n^s,y z_n^t] = [x,y].
\end{align*}
Therefore, central automorphisms act identically on the generators of 
\([\PB_n, \PB_n]\).
\end{proof}

\begin{corollary}
The group {\normalfont \({\rm Aut}_c(\PB_n)\)} preserves both {\normalfont 
\(\Z_n\)} and {\normalfont \(\Bd_n\)}.
\end{corollary}
\begin{proof}
By (\ref{BrunnIsInCommutatorSubgroup}), we have an inclusion \(\Br_n \subseteq 
[\PB_n, \PB_n]\), so \({\rm Aut}_c(\PB_n)\) acts identically on~\(\Z_n\) 
and~\(\Bd_n\).
\end{proof}

Secondly, we prove that \(\Bd_n\) is invariant under \(w_n\). For this, we 
compute the action of~\(w_n\) on \(\{A_{0,j} \mid 1 \le j \le n\}\).

\begin{lemma}
We have
\begin{align}\label{WnAction2}
w_n(A_{0,j}) =
\begin{cases}
A_{j,n}, & 1 \le j \le n-1; \\
A_{0,n} z_n^2, & j = n.
\end{cases}
\end{align}
\end{lemma}
\begin{proof}
For any \(j \in \{1,2,\ldots,n-1\}\), we have
\begin{align*}
w_n(A_{0,j}) &= 
w_n((A_{1,j} A_{2,j} \ldots A_{n-1,j} \cdot A_{n,j})^{-1}) \\ 
&= \left(w_n(A_{1,j}) w_n(A_{2,j}) \ldots w_n(A_{n-1,j}) \cdot w_n(A_{n,j}) 
\right)^{-1} \\
&= \left(A_{1,j} A_{2,j} \ldots A_{n-1,j} \cdot A_{j,n} A_{0,j} A_{j,n}^{-1} 
\right)^{-1} = A_{j,n}.
\end{align*}
It remains to compute \(w_n(A_{0,n})\). For this, we calculate \(w_n(z_n)\).

Since any automorphism preserves the center of the group \(\PB_n\), there 
exists \(s \in \mathbb{Z}\) such that~\(w_n(z_n) = z_n^s\). We claim 
that~\(s=-1\).

To prove the claim, note that the image of the braid \(z_n^s\) under the 
abelianization homomorphism of \(\PB_n\) is
\begin{align*}
A_{1,2}^s (A_{1,3}^s A_{2,3}^s) \ldots (A_{1,n-1}^s \ldots A_{n-2,n-1}^s) 
(A_{1,n}^s \ldots A_{n-1,n}^s).
\end{align*}
Next, we have
\begin{align*}
w_n(z_n) = w_n\left( \prod_{j=2}^n \prod_{i = 1}^{j-1} A_{i,j} \right) &= 
\prod_{j=2}^{n-1} \prod_{i = 1}^{j-1} w_n(A_{i,j}) \cdot \prod_{i = 1}^{n-1} 
w_n(A_{i,n}) \\ &= \prod_{j=2}^{n-1} \prod_{i = 1}^{j-1} A_{i,j} \cdot \prod_{i 
= 1}^{n-1} A_{i,n}A_{0,i}A_{i,n}^{-1}.
\end{align*}
Since the image of \(A_{0,1} A_{0,2} \ldots A_{0,n-1}\) under the 
abelianization homomorphism of \(\PB_n\) is
\begin{align*}
A_{1,2}^{-2} (A_{1,3}^{-2} A_{2,3}^{-2}) \ldots (A_{1,n-1}^{-2} \ldots 
A_{n-2,n-1}^{-2}) (A_{1,n}^{-1} \ldots A_{n-1,n}^{-1}),
\end{align*}
the image of \(w_n(z_n)\) is
\begin{align*}
\prod_{j=2}^{n-1} \prod_{i = 1}^{j-1} A_{i,j}^{1-2} \cdot \prod_{i = 1}^{n-1} 
A_{i,n}^{-1}.
\end{align*}
Therefore, the claim is proved.

Next, we claim that the following identity holds:
\begin{align}\label{ZnIdentity}
A_{0,1} A_{0,2} A_{0,3} \ldots A_{0,n} = z_n^{-2}.
\end{align}
Firstly, the braid~\(\beta \defeq A_{0,1} A_{0,2} A_{0,3} \ldots A_{0,n}\) 
commutes with all the generators~\(\sigma_i\) of~\(\B_n\).
Indeed, Figure~\ref{Ribbons} shows (for~\(n=5\)) that each pair of two 
adjacent strands of~\(\beta\) bounds a ribbon, and it suffices to push the 
additional crossing~\(\sigma_i\) through the ribbon to establish the 
equality~\(\sigma_i \beta = \beta \sigma_i\). Next, recall that the center 
of~\(\B_n\) agrees with those of~\(\PB_n\) (see~\cite[Theorem~15]{GJP15}).
Thus, there exists~\(t \in \mathbb{Z}\) such that~\(\beta = z_n^t\). The 
images of~\(\beta\) and~\(z_n^{-2}\) under the abelianization homomorphism are 
equal, because each pair of two adjacent strands of~\(\beta\), as well as 
of~\(z_n^{-2}\), has precisely four negative half-twists. Therefore,~\(t=-2\).

\begin{figure}
\centering
\includegraphics[height = 5cm]{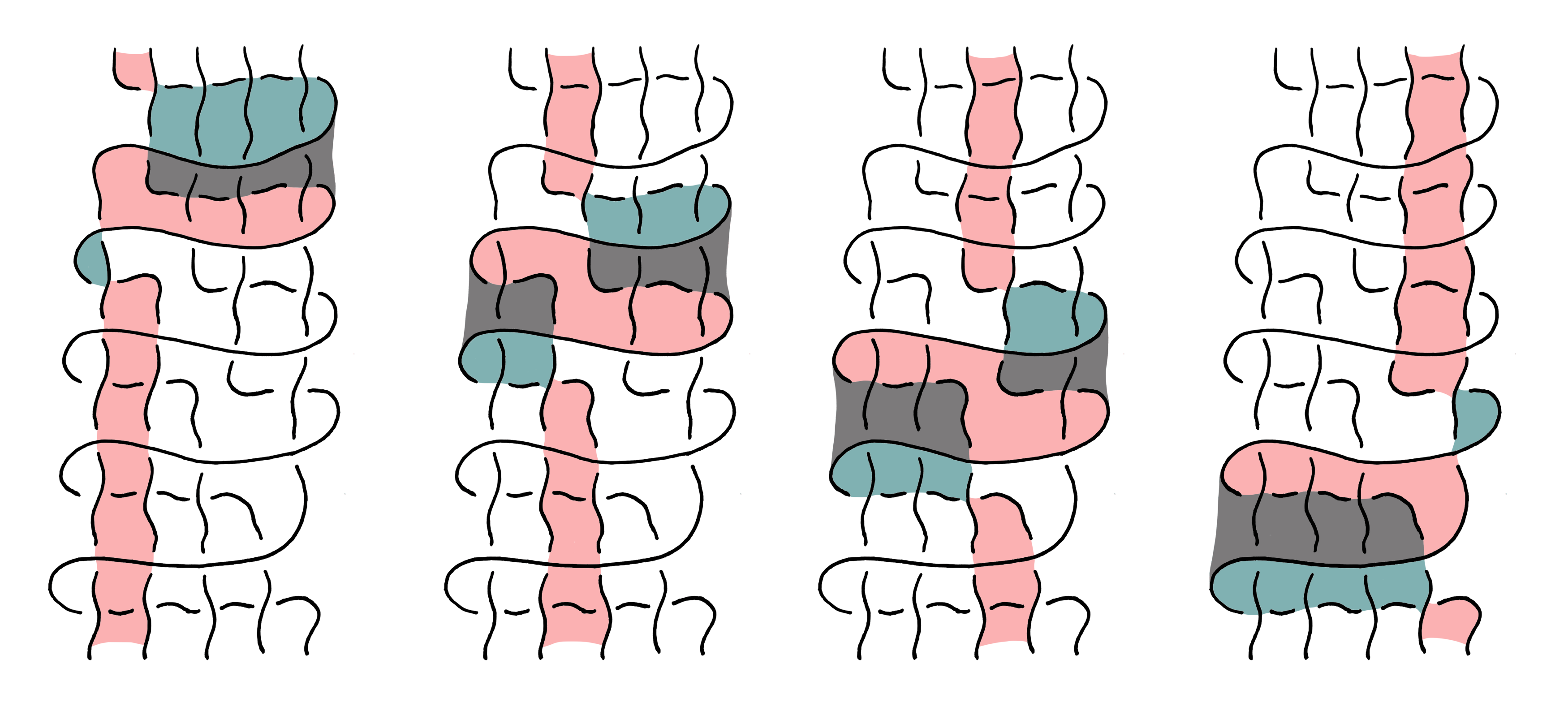}
\caption{Pairs of adjacent strands of \(A_{0,1} A_{0,2} \ldots A_{0,n}\) bound 
ribbons.}
\label{Ribbons}
\end{figure}

Finally, we have
\[w_n(A_{0,n}) = w_n(A_{0,n-1}^{-1} \ldots A_{0,1}^{-1} z_n^{-2}) = 
A_{n-1,n}^{-1} \ldots A_{1,n}^{-1} z_n^2 = A_{0,n} z_n^2. \qedhere\]
\end{proof}
\medskip

\begin{corollary}
The automorphism {\normalfont \(w_n\)} preserves {\normalfont \(\Bd_n\)}.
\end{corollary}
\begin{proof}
We have
\begin{align*}
w_n\left(\Bd_n\right) 
\overset{(\ref{BoundaryBrunnianSymmetricCommutator})}&{=} w_n\left([\llangle 
A_{0,1} \rrangle, \llangle A_{0,2} \rrangle, \ldots, \llangle A_{0,n-1} 
\rrangle, \llangle A_{0,n} \rrangle]_S\right) \\
&= [\llangle w_n(A_{0,1}) \rrangle, \llangle w_n(A_{0,2}) \rrangle, \ldots, 
\llangle w_n(A_{0,n-1}) \rrangle, \llangle w_n(A_{0,n}) \rrangle]_S \\
\overset{(\ref{WnAction2})}&{=} [\llangle A_{1,n} \rrangle, \llangle A_{2,n} 
\rrangle, \ldots, \llangle A_{n-1,n} \rrangle, \llangle A_{0,n} z_n^2 
\rrangle]_S \\
&= [\llangle A_{1,n} \rrangle, \llangle A_{2,n} \rrangle, \ldots, \llangle 
A_{n-1,n} \rrangle, \llangle A_{0,n} \rrangle]_S \\
&= [\llangle A_{0,n} \rrangle, \llangle A_{1,n} \rrangle, \llangle A_{2,n} 
\rrangle, \ldots, \llangle A_{n-1,n} \rrangle]_S \\
&=\partial_{n+1}\left([\llangle A_{1,n+1} \rrangle, \llangle A_{2,n+1} 
\rrangle, \ldots, \llangle A_{n,n+1} \rrangle]_S\right) \\
\overset{(\ref{BrunnianSymmetricCommutator})}&{=} 
\partial_{n+1}\left(\Br_{n+1}\right) = \Bd_n. \qedhere
\end{align*}
\end{proof}
\medskip

Thirdly, we prove that \(\Z_n\) is invariant under \(w_n\). For this, we give 
a generating set for the group \({\rm Ker}(\partial_n)\).

\begin{lemma}\label{Lemma1}
The set {\normalfont \(\{A_{0,j} \mid 2 \le j \le n\}\)} generates 
{\normalfont \({\rm Ker}(\partial_n)\)}.
\end{lemma}
\begin{proof}
For \(j \in \{2,3,\ldots,n\}\), we have
\begin{align*}
\theta_n(A_{0,j}) &= 
\theta_n((A_{1,j} \cdot A_{2,j} \ldots A_{n,j})^{-1}) \\
&= \left(\theta_n(A_{1,j}) \cdot \theta_n(A_{2,j}) \ldots \theta_n(A_{n,j}) 
\right)^{-1} \\
&=\left(A_{1,j}^{-1} A_{0,j} A_{1,j} \cdot A_{2,j} \ldots A_{n,j} \right)^{-1} 
= A_{1,j}.
\end{align*}
Next,
\begin{align*}
{\rm Ker}(\partial_n) &= {\rm Ker}(d_1 \circ \theta_n) = \theta_n^{-1}({\rm 
Ker}(d_1)) \\
\overset{(\ref{DkGenerators})}&{=} \theta_n^{-1}(\langle A_{1,j} \mid 2 \le j 
\le n\rangle) = \langle \theta_n^{-1}(A_{1,j}) \mid 2 \le j \le n\rangle \\
&=\langle A_{0,j} \mid 2 \le j \le n \rangle.\qedhere
\end{align*}
\end{proof}
\medskip

\begin{corollary}
The automorphism {\normalfont \(w_n\)} preserves {\normalfont \(\Z_n\)}.
\end{corollary}
\begin{proof}
Firstly, we prove that 
\begin{align*}
w_n(\Z_n) \subseteq \bigcap_{k=1}^{n-1} {\rm Ker}(d_k). 
\end{align*}
Let \(k \in \{1,2,\ldots,n-1\}\). For any \(j \in \{1,2,\ldots,n-1\}\), we 
have \[w_n(A_{k,j}) = A_{k,j} \in {\rm Ker}(d_k).\] Besides, \[w_n(A_{k,n}) = 
A_{k,n}A_{0,k}A_{k,n}^{-1} \in {\rm Ker}(d_k).\] It follows that~\(w_n({\rm 
Ker}(d_k)) \subseteq {\rm Ker}(d_k)\).
Therefore,
\begin{align*}
w_n(\Z_n) \subseteq w_n\left(\bigcap_{k=1}^{n-1} {\rm Ker}(d_k)\right) 
\subseteq \bigcap_{k=1}^{n-1} {\rm Ker}(d_k).
\end{align*}

Secondly, we prove that~\(w_n(\Z_n) \subseteq {\rm Ker}(d_n)\). 
Let~\(\beta \in \Z_n\). 
By Lemma~\ref{Lemma1}, there exists a word~\(u_1\) in the 
alphabet~\(\{A_{0,2}^{\pm 1}, A_{0,3}^{\pm 1}, \ldots, A_{0,n-1}^{\pm 1}, 
A_{0,n}^{\pm 1}\}\) that represents the braid~\(\beta\).
We claim that the exponent sum of the letter~\(A_{0,n}\) in~\(u_1\) is zero.

Indeed,~\(\theta_n(\beta)\) can be presented as a word~\(u_2\) in the alphabet
\[\{A_{1,2}^{\pm 1}, \ldots, A_{1,n-1}^{\pm 1}, A_{1,n}^{\pm 1}\}.\]
Recall that~\(\Z_n \subseteq [\PB_n,\PB_n]\). 
Hence,~\(\theta_n(\beta) \in [\PB_n,\PB_n]\). 
It follows that the exponent sum of the letter~\(A_{1,n}\) in the word~\(u_2\) 
is equal to zero.
Thus the exponent sum of~\(A_{0,n}\) in~\(u_1\) is zero. The claim is proved.

By (\ref{WnAction2}), the braid \(w_n(\beta)\) can be presented as a word 
\(u_3\) in the alphabet
\[\{A_{2,n}^{\pm 1}, A_{3,n}^{\pm 1}, \ldots, A_{n-1,n}^{\pm 1}, L^{\pm 1}\},\]
where \(L \defeq A_{0,n} z_n^2\). By the above claim, the exponent sum of 
\(L\) in~\(u_3\) is zero. Since the braid \(z_n\) is central, \(w_n(\beta)\) 
can be presented as a word~\(u_4\) in the alphabet
\[\{A_{2,n}^{\pm 1}, A_{3,n}^{\pm 1}, \ldots, A_{n-1,n}^{\pm 1}, A_{0,n}^{\pm 
1}\}.\]
Note that for any \(k \in \{2,3,\ldots,n-1,0\}\), we have~\(A_{k,n} \in {\rm 
Ker}(d_n)\). It follows that \(w_n(\beta) \in {\rm Ker}(d_n)\).

Thirdly, we prove that \(w_n(\Z_n) \subseteq {\rm Ker}(\partial_n)\). We have
\begin{align*}
w_n(\Br_n) \overset{(\ref{BrunnianSymmetricCommutator})}&{=} 
w_n\left([\llangle A_{1,n} \rrangle, \llangle A_{2,n} \rrangle, \ldots, 
\llangle A_{n-1,n} \rrangle]_S\right) \\
&= [\llangle w_n(A_{1,n}) \rrangle, \llangle w_n(A_{2,n}) \rrangle, \ldots, 
\llangle w_n(A_{n-1,n}) \rrangle]_S \\
&= [\llangle A_{0,1} \rrangle, \llangle A_{0,2} \rrangle, \ldots, \llangle 
A_{0,n-1} \rrangle]_S.
\end{align*}
By part (1) of~\cite[Lemma 6.5.2]{BCWW06}, for any~\(j \in \{2,3,\ldots,n\}\), 
one has \(\partial_n(A_{0,j}) = \triv\). Hence,
\begin{align*}
\partial_n\left(w_n(\Br_n)\right) &= \partial_n\left( [\llangle A_{0,1} 
\rrangle, \llangle A_{0,2} \rrangle, \ldots, \llangle A_{0,n-1} \rrangle]_S 
\right) \\
&= [\llangle \partial_n(A_{0,1}) \rrangle, \llangle \partial_n(A_{0,2}) 
\rrangle, \ldots, \llangle \partial_n(A_{0,n-1}) \rrangle]_S \\
&= [\llangle \partial_n(A_{0,1}) \rrangle, \llangle \triv \rrangle, \ldots, 
\llangle \triv \rrangle]_S = \{\triv\}.
\end{align*}
Therefore, \(w_n(\Z_n) \subseteq w_n(\Br_n) \subseteq {\rm Ker}(\partial_n)\).
\end{proof}

Fourthly, we prove that both~\(\Z_n\) and~\(\Bd_n\) are invariant under~\({\rm 
Aut}(\B_n)\).

Recall that~\({\rm Aut}(\B_n)\) is generated by the reflection 
automorphism~\(\chi_n\) and the conjugations by all elements of~\(\B_n\). As 
indicated in~\cite[Lemma 4.3]{LW09}, the invariance under~\(\chi_n\) would 
follow from the fact that~\(\chi_n\) commutes with the 
homomorphisms~\(d_1,\ldots,d_n,\partial_n\colon \PB_n \to \PB_{n-1}\). The 
identities~\(d_k \circ \chi_n = \chi_{n-1} \circ d_k\) are obvious. To prove 
the identity~\(\partial_n \circ \chi_n = \chi_{n-1} \circ \partial_n\), 
in~\cite[Lemma 2.3]{LW09} Li and Wu suggest drawing pictures. For the sake of 
completeness, we follow their suggestion and provide proof.

\begin{figure}
\centering
\includegraphics[height = 5.5cm]{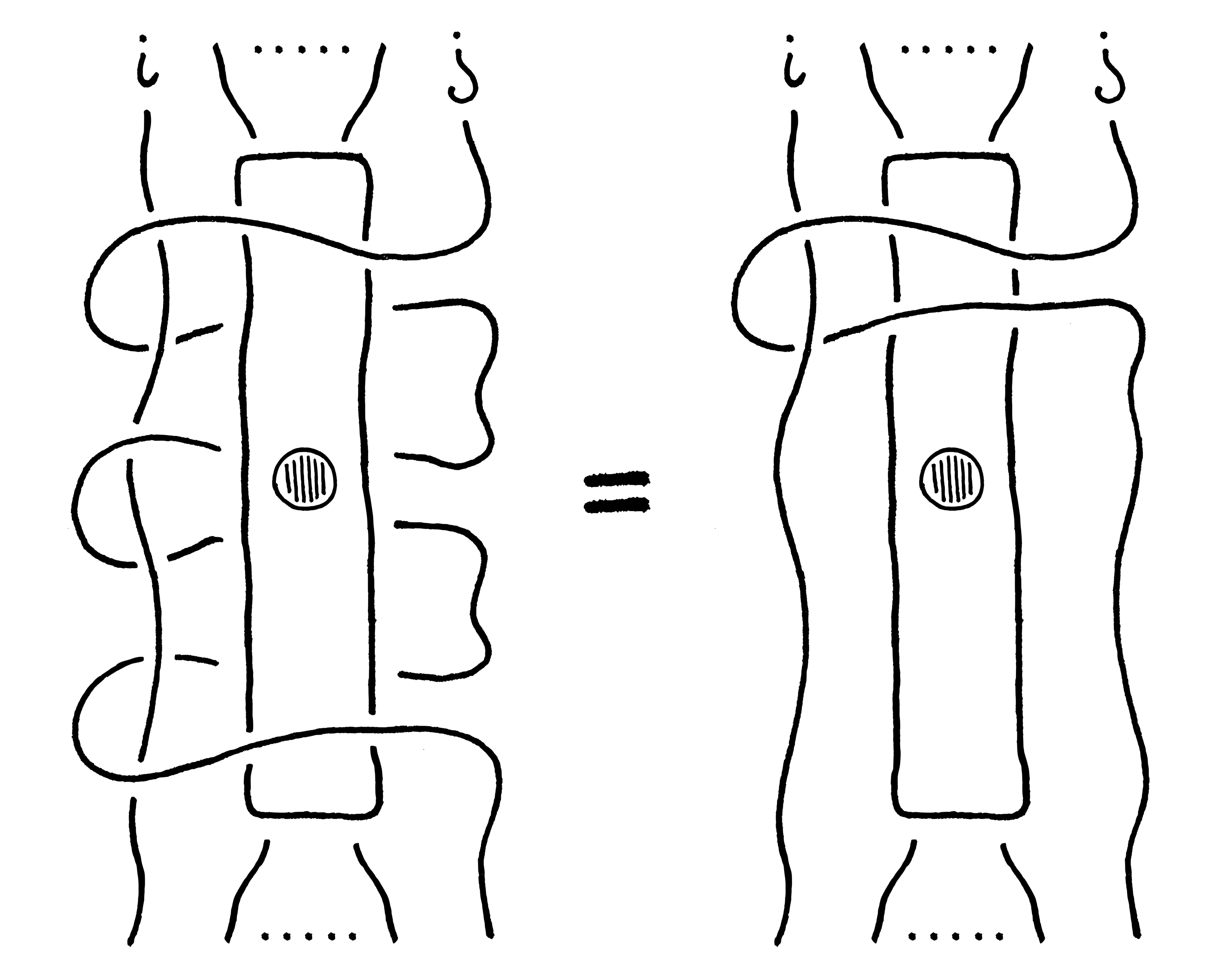}
\caption{Conjugation between~\(A_{i,j}^{-1}\) and~\(\chi_n(A_{i,j})\).}
\label{ConjugationAij}
\end{figure}

\begin{figure}[ht]
\centering
\includegraphics[height = 5cm]{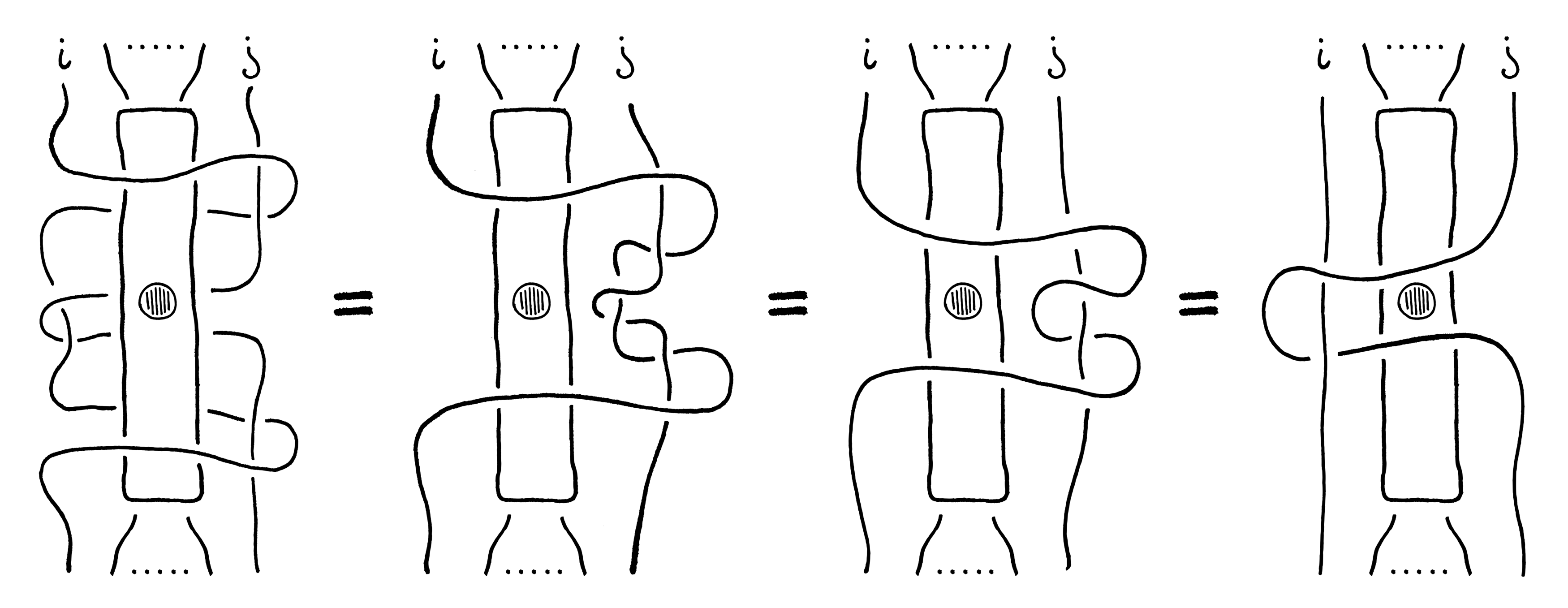}
\caption{Alternative conjugation between~\(A_{i,j}^{-1}\) 
and~\(\chi_n(A_{i,j})\).}
\label{AltConjugationAij}
\end{figure}

\begin{figure}
\centering
\includegraphics[height = 5cm]{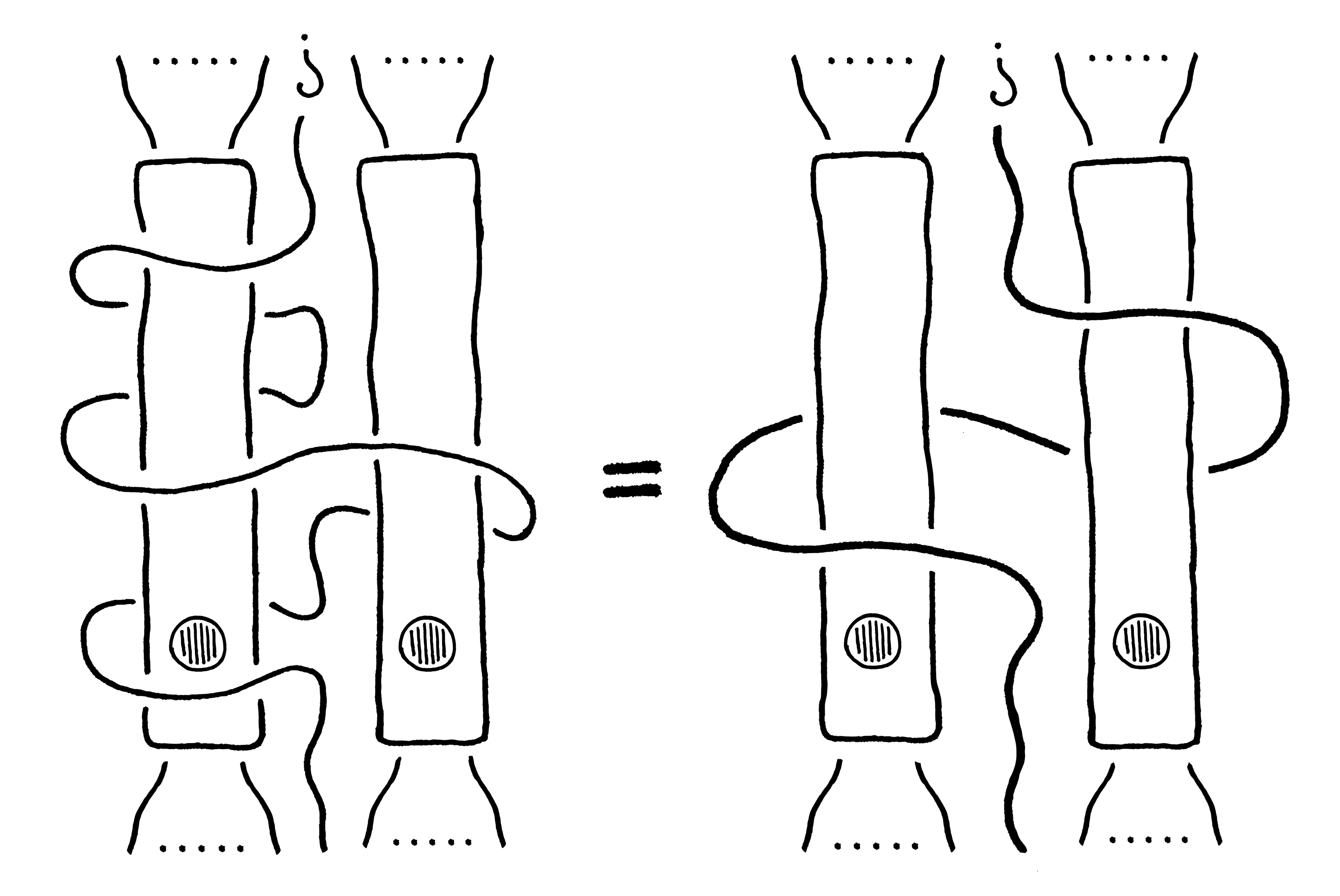}
\caption{Conjugation between~\(A_{0,j}^{-1}\) and~\(\chi_n(A_{0,j})\).}
\label{ConjugationA0j}
\end{figure}

Below, we express braids of the form \(\chi_n(A_{i, j})\) in terms of the 
standard generators of~\(\PB_n\). See Figures \ref{ConjugationAij}, 
\ref{AltConjugationAij}, and~\ref{ConjugationA0j} for a geometric 
interpretation. Here, the~\(j\)th strand of the braid~\((A_{i,j}A_{i+1,j} 
\ldots A_{j-1,j})^{-1}\) is wrapped around the strands 
\(i,{i+1},\ldots,{j-1}\), and the~\(i\)th strand of the braid~\((A_{i,i+1} 
A_{i,i+2} \ldots A_{i,j})\) is wrapped around the strands 
\({i+1},{i+2},\ldots,j\).

\begin{lemma}\label{ReflecionExplicitFormulas}
For any \(i, j \in \{1, 2, \ldots, n\}\) such that \(i < j\), one has
\begin{align}
\chi_n(A_{i,j}) &= (A_{i,j} A_{i+1,j} \ldots A_{j-1,j})^{-1} A_{i,j}^{-1} 
(A_{i,j} A_{i+1,j} \ldots A_{j-1,j}), \label{ChiEq1} \\
\chi_n(A_{i,j}) &= (A_{i,i+1} A_{i,i+2} \ldots A_{i,j}) A_{i,j}^{-1} 
(A_{i,i+1} A_{i,i+2} \ldots A_{i,j})^{-1}, \label{ChiEq2}
\end{align}
and for any \(j \in \{1,2,\ldots,n\}\), one has
\begin{equation}\label{ChiEq3}
\chi_n(A_{0,j}) = (A_{1,j} A_{2,j} \ldots A_{j-1,j})^{-1} A_{0,j}^{-1} 
(A_{1,j} A_{2,j} \ldots A_{j-1,j}).
\end{equation}
\end{lemma}

\begin{proof}
Equality (\ref{ChiEq1}) is satisfied by \cite[Lemma 3.1]{BNS18}.
It is easy to see that
\begin{align}\label{AijAlternative}
A_{s,t} = \sigma_s^{-1} \sigma_{s+1}^{-1} \ldots \sigma_{t-2}^{-1} 
\sigma_{t-1}^2 \sigma_{t-2} \ldots \sigma_{s+1} \sigma_s
\end{align}
for any \(s, t \in \{1,2,\ldots,n\}\) such that \(s < t\).
Let \(\Phi\) be the involutive automorphism of \(\B_n\) given by 
\(\Phi(\sigma_m) = \sigma_{n-m}^{-1}\). Then
\begin{align*}
\Phi(A_{s,t}) &= \Phi(\sigma_{t-1} \ldots \sigma_{s+1} \sigma_s^2 
\sigma_{s+1}^{-1} \ldots \sigma_{t-1}^{-1}) \\
&= \sigma_{n-t+1}^{-1} \ldots \sigma_{n-s-1}^{-1} \sigma_{n-s}^{-2} 
\sigma_{n-s-1} \ldots \sigma_{n-t+1} \\
&= (\sigma_{n-t+1}^{-1} \ldots \sigma_{n-s-1}^{-1} \sigma_{n-s}^2 
\sigma_{n-s-1} \ldots \sigma_{n-t+1})^{-1} \overset{(\ref{AijAlternative})}{=} 
A_{n-t+1,n-s+1}^{-1}.
\end{align*}
Let \(i' \defeq n - j + 1\) and \(j' \defeq n - i + 1\). Since~\(\Phi \circ 
\chi_n = \chi_n \circ \Phi\), one has
\begin{align*}
\chi_n(A_{i,j}) &= \chi_n(\Phi(A_{i',j'}^{-1})) = 
\Phi(\chi_n(A_{i',j'}^{-1})) \\
&\overset{(\ref{ChiEq1})}{=} \Phi\left((A_{i',j'} A_{i'+1,j'} \ldots 
A_{j'-1,j'})^{-1} A_{i',j'} (A_{i',j'} A_{i'+1,j'} \ldots A_{j'-1,j'})\right) \\
&= (A_{i,j}^{-1} A_{i,j-1} 
\ldots A_{i,i+1}^{-1})^{-1} A_{i,j}^{-1} (A_{i,j}^{-1} A_{i,j-1} 
\ldots A_{i,i+1}^{-1}) \\
&= (A_{i,i+1} A_{i,i+2} \ldots A_{i,j}) A_{i,j}^{-1} (A_{i,i+1} A_{i,i+2} 
\ldots A_{i,j})^{-1},
\end{align*}
hence the identity (\ref{ChiEq2}) follows.

To show that the identity (\ref{ChiEq3}) holds, we generously apply the 
results just proved:
\begin{align*}
\chi_n(A_{0,j}) &= \chi_n(\underbrace{A_{1,j} A_{2,j} \ldots 
A_{j-1,j}}_{\text{apply (\ref{ChiEq1})}} \cdot \underbrace{A_{j,j+1} A_{j,j+2} 
\ldots A_{j,n}}_{\text{apply (\ref{ChiEq2})}})^{-1} \\
&= \left[\left(A_{1,j} A_{2,j} \ldots A_{j-1,j}\right)^{-1}\cdot 
\left.\highlight[MyRed]{A_{1,j}}\right.^{-1}\cdot 
\left.\highlight[MyRed]{A_{1,j}} \highlight[MyGreen]{A_{2,j} \ldots 
A_{j-1,j}}\right.\right. \\
& \cdot\left.\highlight[MyGreen]{A_{2,j} A_{3,j} \ldots 
A_{j-1,j}}\right.^{-1}\cdot \left.\highlight[MyRed]{A_{2,j}}\right.^{-1}\cdot 
\left.\highlight[MyRed]{A_{2,j}} \highlight[MyGreen]{A_{3,j} \ldots 
A_{j-1,j}}\right.\cdot \ldots \\
& \cdot \left.\highlight[MyGreen]{A_{j-1,j}}\right.^{-1}\cdot 
\left.\highlight[MyRed]{A_{j-1,j}}\right.^{-1}\cdot 
\left.\highlight[MyRed]{A_{j-1,j}}\right. \cdot \highlight[MyGreen]{A_{j,j+1}} 
\cdot \left.\highlight[MyGreen]{A_{j,j+1}}\right.^{-1}\cdot 
\left.\highlight[MyRed]{A_{j,j+1}}\right.^{-1} \\
& \cdot \left.\highlight[MyRed]{A_{j,j+1}} 
\highlight[MyGreen]{A_{j,j+2}}\right.\cdot 
\left.\highlight[MyGreen]{A_{j,j+2}}\right.^{-1}\cdot 
\left.\highlight[MyRed]{A_{j,j+1} A_{j,j+2}}\right.^{-1}\cdot \ldots \\
& \left.\cdot\left.\highlight[MyRed]{A_{j,j+1} A_{j,j+2} \ldots A_{j,n-1}} 
\highlight[MyGreen]{A_{j,n}}\right. \cdot 
\left.\highlight[MyGreen]{A_{j,n}}\right.^{-1}\cdot \left(A_{j,j+1} A_{j,j+2} 
\ldots A_{j,n}\right)^{-1}\right]^{-1} \\
&= A_{j,j+1} A_{j,j+2} \ldots A_{j,n} \cdot A_{1,j} A_{2,j} \ldots A_{j-1,j} \\
&= (A_{1,j} \ldots A_{j-1,j})^{-1} (A_{1, j} \ldots A_{j-1,j} \cdot A_{j, j+1} 
\ldots A_{j,n}) (A_{1,j} \ldots A_{j-1,j}) \\
&= (A_{1,j} \ldots A_{j-1,j})^{-1} A_{0,j}^{-1} (A_{1,j} \ldots A_{j-1,j}).
\end{align*}
Here adjacent same-colored blocks cancel each other.
\end{proof}

\begin{lemma}\label{ReflectionCommutesWithPartial}
The following equality holds:
\begin{align*}
\partial_n \circ \chi_n = \chi_{n-1} \circ \partial_n.
\end{align*}
\end{lemma}
\begin{proof}
Recall that~\(\partial_n = d_1 \circ \theta_n\). %

If~\(i \ge 2\), then~\(\theta_n(A_{i,j}) = A_{i,j}\), and by~\eqref{ChiEq1}, 
we have~\(\theta_n(\chi_n(A_{i,j})) = \chi_n(A_{i,j})\). Consequently,
\begin{align*}
(\theta_n \circ \chi_n)(A_{i,j}) = \chi_n(A_{i,j}) = (\chi_n \circ 
\theta_n)(A_{i,j}).
\end{align*}
Therefore, 
\begin{align*}
(\partial_n \circ \chi_n)(A_{i,j}) &= (d_1 \circ \theta_n \circ 
\chi_n)(A_{i,j}) = (d_1 \circ \chi_n \circ \theta_n)(A_{i,j}) \\
 &= (\chi_{n-1} \circ d_1 \circ \theta_n)(A_{i,j}) = (\chi_{n-1} \circ 
\partial_n)(A_{i,j}).
\end{align*}

Next, for~\(i = 1\), we have
\begin{align*}
\chi_n(A_{1,j}) \overset{(\ref{ChiEq1})}{=} &(A_{1,j} A_{2,j} \ldots 
A_{j-1,j})^{-1} A_{1,j}^{-1} (A_{1,j} A_{2,j} \ldots A_{j-1,j}) \\
 =\ &(A_{2,j} A_{3,j} \ldots A_{j-1,j})^{-1} A_{1,j}^{-1} (A_{2,j} A_{3,j} 
\ldots A_{j-1,j}).
\end{align*}
Recall that \(\theta_n(A_{1,j}) = A_{1,j}^{-1}A_{0,j}A_{1,j}\). We have
\begin{align*}
\theta_n(\chi_n(A_{1,j})) &= \theta_n\left((A_{2,j} A_{3,j} \ldots 
A_{j-1,j})^{-1} A_{1,j}^{-1} (A_{2,j} A_{3,j} \ldots A_{j-1,j})\right) \\
 &= (A_{1,j} A_{2,j} \ldots A_{j-1,j})^{-1} A_{0,j}^{-1} (A_{1,j} A_{2,j} 
\ldots A_{j-1,j}) \\
&\overset{(\ref{ChiEq3})}{=} \chi_n(A_{0,j}).
\end{align*}
Therefore, 
\begin{align*}
(\partial_n \circ \chi_n)(A_{1,j}) &= d_1(\chi_n(A_{0,j})) = 
\chi_{n-1}(d_1(A_{0,j})) \\
 &= \chi_{n-1}(d_1(A_{1,j}^{-1}A_{0,j}A_{1,j})) = 
\chi_{n-1}(d_1(\theta_n(A_{1,j}))) \\
 &= (\chi_{n-1} \circ \partial_n)(A_{1,j}).\qedhere
\end{align*}
\end{proof}
\medskip

Finally, to prove that the group \({\rm Aut}(\B_n)\) preserves both \(\Z_n\) 
and \(\Bd_n\), it remains to show that these subgroups are invariant under the 
conjugations by all elements of the braid group \(\B_n\). The subgroup 
\(\Bd_n\) of \(\B_n\) is normal by \cite[Lemma 3.10]{LW09}.

\begin{lemma}
The subgroup {\normalfont \(\Z_n\)} of {\normalfont \(\B_n\)} is normal.
\end{lemma}
\begin{proof}
Given \(\beta \in \Z_n\) and \(\alpha \in \B_n\), we aim to show that 
\(\Psi_\alpha(\beta) \in \Z_n\).

By \eqref{PermutationOfDs}, the automorphism \(\Psi_\alpha\) permutes the set 
\(\{{\rm Ker}(d_1), {\rm Ker}(d_1), \ldots, {\rm Ker}(d_n)\}\).
Hence, \(\Psi_\alpha(\beta) \in \Br_n\).

We aim to prove that \(\Psi_\alpha(\beta) \in {\rm Ker}(\partial_n)\).
Without loss of generality, \(\alpha\) is the element of the generating set 
\(\{\sigma_1^{\pm 1}, \ldots, \sigma_{n-1}^{\pm 1}\}\) of \(\B_n\).
Since \(\Z_n\) is invariant under the reflection automorphism \(\chi_n\) of 
\(\B_n\) and \(\Psi_{\sigma_k} = \chi_n \circ \Psi_{\sigma_k^{-1}} \circ 
\chi_n\),
we can assume that \(\alpha = \sigma_k^{-1}\), where \(k \in 
\{1,2,\ldots,n-1\}\).

By Lemma \ref{Lemma1}, there exists a word~\(u_1\) in the 
alphabet~\(\{A_{0,2}^{\pm 1}, A_{0,3}^{\pm 1}, \ldots, A_{0,n-1}^{\pm 1}, 
A_{0,n}^{\pm 1}\}\) that represents the braid \(\beta\). By \cite[Lemma 
3.1]{LW09}, we have
\begin{align*}
\Psi_{\sigma_k^{-1}}(A_{0,j}) =
\begin{cases}
A_{0,j}, & k \in \{1,2,\ldots,j-2, j+1, \ldots, n\}; \\
A_{0,j}^{-1}A_{0,j-1}A_{0,j}, & k = j-1; \\
A_{0,j+1}, & k = j.
\end{cases}
\end{align*}
By part (1) of~\cite[Lemma 6.5.2]{BCWW06}, for any \(j \in \{2,3,\ldots,n\}\), 
one has \(\partial_n(A_{0,j}) = \triv\). In particular,
for any \(k \in \{2,3,\ldots,n-1\}\) and \(j \in \{2,3,\ldots,n\}\), we have 
\begin{align*}
\partial_n\left(\Psi_{\sigma_k^{-1}}(A_{0,j}) \right) = \triv.
\end{align*}
Therefore, it is enough to prove the statement for \(k = 1\).

Let \(m\) be the exponent sum of the letter \(A_{0,2}\) in \(u_1\). Firstly, 
we have
\begin{align}\label{PartialA01}
\partial_n(A_{0,1}) = z_{n-1}^2.
\end{align}
Indeed,
\begin{align*}
\partial_n(A_{0,1}) &= \left( d_1 \circ \theta_n \right) \left( A_{1,n}^{-1} 
A_{1,n-1}^{-1} \ldots A_{1,2}^{-1} \right) \\
&= d_1 \left( \left(A_{1,n}^{-1}A_{0,n} A_{1,n}\right)^{-1} \left( 
A_{1,n-1}^{-1}A_{0,n-1} A_{1,n-1}^{-1} \right)^{-1} \ldots \left( 
A_{1,2}^{-1}A_{0,2} A_{1,2} \right)^{-1}\right) \\
&= A_{0,n-1}^{-1} A_{0,n-2}^{-1} \ldots A_{0,1}^{-1} 
\overset{(\ref{ZnIdentity})}{=} z_{n-1}^2.
\end{align*}
Secondly, we have 
\begin{align*}
\partial_n(\Psi_{\sigma_1^{-1}}(\beta)) = 
\partial_n(\Psi_{\sigma_1^{-1}}(A_{0,2}))^m = \partial_n(A_{0,1})^m = 
z_{n-1}^{2m}.
\end{align*}
Next, we claim that
\(m=0\).

To prove the claim, note that the braid \(\theta_n(\beta)\) can be presented 
as a word \(u_2\) in the alphabet~\(\{A_{1,2}^{\pm 1}, A_{1,3}^{\pm 1}, \ldots, 
A_{1,n-1}^{\pm 1}, A_{1,n}^{\pm 1}\}\). Recall that \(\Z_n \subseteq 
[\PB_n,\PB_n]\). In particular, we have~\(\theta_n(\beta) \in [\PB_n,\PB_n]\). 
It follows that the exponent sum of the letter~\(A_{1,2}\) in the word~\(u_2\) 
is equal to zero. Thus the exponent sum of~\(A_{0,2}\) in~\(u_1\) is zero. The 
claim is proved.
\end{proof}

\section{Proof of Proposition \ref{Proposition1}}

Firstly, we show that \(\Br_3 = \Z_3 = [\PB_3,\PB_3]\).
It is clear that~\([\PB_3, \PB_3] \subseteq \Br_3\), and 
by~(\ref{BrunnIsInCommutatorSubgroup}), we have \([\PB_3, \PB_3] = \Br_3\).
Next,~\([\PB_3, \PB_3] \subseteq \Z_3\) because~\(\theta_3([\PB_3, \PB_3]) = 
[\PB_3, \PB_3]\) and~\(d_1([\PB_3, \PB_3]) = \{\triv\}\).
Therefore, \(\Z_3 = [\PB_3, \PB_3]\).

Secondly, we find a basis of the free group \(\Z_3\).

\begin{proposition}
The set {\normalfont \(\{[A_{1,2}^n, A_{2,3}^m] \mid n,m \in \mathbb{Z} 
\backslash \{0\}\}\)} is a basis of the free group~{\normalfont \(\Z_3\)}.
\end{proposition}
\begin{proof}
By \cite[Section~3.5]{GJP15},
the subgroup \(F_2 \defeq \langle A_{1,2}, A_{2,3} \rangle\) is free.
By \cite[p.~196, ex.~24]{MKS76}, the set \(\{[A_{1,2}^n, A_{2,3}^m] \mid n,m 
\in \mathbb{Z} \backslash \{0\}\}\) is a basis of \([F_2,F_2]\). To prove the 
proposition, it is enough to show that~\([\PB_3, \PB_3] = [F_2, F_2]\).

The inclusion \([\PB_3, \PB_3] \supseteq [F_2, F_2]\) is obvious. Next, 
\[[\PB_3, \PB_3] = \Br_3 \subseteq {\rm Ker}(d_2) = F_2.\]
Recall that the quotient group~\(F_2 \slash [F_2, F_2]\) is free abelian of 
rank~\(2\) and the canonical projection \(F_2 \to F_2 \slash [F_2, F_2]\) 
assigns to an element of~\(F_2\), given by a word in the 
alphabet~\(\{A_{1,2}^{\pm 1}, A_{2,3}^{\pm 1}\}\), the pair of exponent sums of 
the letters \(A_{1,2}\) and \(A_{2,3}\).
By definition, for any element of \([\PB_3, \PB_3]\),
the exponent sums of the letters \(A_{1,2}\) and \(A_{2,3}\) are zero. 
Therefore,~\([\PB_3, \PB_3] \subseteq [F_2,F_2]\).
\end{proof}

Thirdly, we find a generator of the quotient group \(\Z_3 \slash \Bd_3\).

Let \(a \defeq A_{1,2}\), \(b \defeq A_{2,3}\), and \(c \defeq A_{1,3}\). We 
need the following identities:
\begin{align}
c^{-1} a c &= b a b^{-1}, \label{Conj1} \\
c a c^{-1} &= a^{-1}b^{-1} a ba, \label{Conj2} \\
c^{-1} b c &= ba b a^{-1} b^{-1}, \label{Conj3} \\
c b c^{-1} &= a^{-1} b a. \label{Conj4}
\end{align}
Denote by \(\pi \colon \Z_3 \to \Z_3\slash \Bd_3\) the canonical projection. 
Given \(x,y \in \Z_3\), we write \(x \equivv y\) if~\(\pi(x) = \pi(y)\).

\begin{lemma}\label{CyclicityOfPi3}
The image of {\normalfont \([a,b]\)} under {\normalfont \(\pi\)} generates the 
quotient group {\normalfont \(\Z_3\slash \Bd_3\)}.
\end{lemma}
\begin{proof}
It is enough to show that for all non-zero \(n,m \in \mathbb{Z}\), one has 
\([a^n,b^m] \equivv [a,b]^{nm}\).

Since \(\Bd_3\) is a characteristic subgroup of \(\PB_3\), for any 
automorphism \(f\) of \(\PB_3\), the equivalence \(x \equivv y\) implies \(f(x) 
\equivv f(y)\). Below, we use the following identities, which are easy to prove:
\begin{align*}
\Psi_a([a^n,b^m]) &=
[a^{n+1}, b^m] [a, b^m]^{-1}, \\
\Psi_{a}^{-1}([a^n,b^m]) &=
[a^{n-1}, b^m] [a^{-1}, b^m]^{-1},\\
\Psi_b([a^n, b^m]) &=
[a^n,b]^{-1} [a^n, b^{m+1}], \\
\Psi_{b}^{-1}([a^n,b^m]) &=
[a^n,b^{-1}]^{-1} [a^n, b^{m-1}].
\end{align*}

Firstly, we claim that for any non-zero \(n \in \mathbb{Z}\), we have \([a^n, 
b^2] \equivv [a^n, b]^2\). The proof is by induction on \(|n|\).

The base case is \(|n|=1\).
By \cite[Proposition 3.15]{LW09}, one has
\[[[A_{1,2}, A_{2,3}], A_{0,1}] \in \Bd_3.\]
Recall that~\(A_{0,1} = c^{-1} a^{-1}\). By (\ref{Conj2}) and (\ref{Conj4}), 
we have
\begin{align*}
[[A_{1,2}, A_{2,3}], A_{0,1}] &= [[a,b], c^{-1}a^{-1}] = [a,b]^{-1} a c [a,b] 
c^{-1}a^{-1} \\
&= [a,b]^{-1} a [cac^{-1},cbc^{-1}] a^{-1} \\
&= [a,b]^{-1} a [a^{-1} b^{-1} a ba, a^{-1} b a] a^{-1} \\
&= [a,b]^{-1} [b^{-1} a b, b] \\
&= [a,b]^{-1} [b,a] [a,b^2] = [b,a]^2 [a,b^2].
\end{align*}
Hence, \([a,b^2] \equivv [a,b]^2\). It follows that
\begin{align*}
[a^{-1}, b^2]^{-1} = \Psi_{a}^{-1}([a,b^2]) \equivv \Psi_a^{-1}([a,b]^2) = 
[a^{-1}, b]^{-2}.
\end{align*}
Therefore, the statement holds for \(|n|=1\).

Next, we prove the inductive step. Since the quotient \(\Z_3 \slash \Bd_3\) is 
abelian (see \cite[Proposition 4.1.3]{BCWW06}), one has \((xy)^m \equivv x^m 
y^m\) for any \(x,y \in \Z_3\) and \(m \in \mathbb{Z}\). For \(n \ge 1\), we 
have
\begin{align*}
[a^{n+1}, b^2] [b, a]^2 &\equivv [a^{n+1}, b^2] [b^2, a] = \Psi_{a}([a^n,b^2]) 
\equivv \Psi_{a}([a^n,b])^2 = \left([a^{n+1}, b] [b, a] \right)^2 \\
&\equivv [a^{n+1}, b]^2 [b, a]^2.
\end{align*}
The above shows that the statement holds for all \(n \ge 1\). Next, for \(n 
\le -1\), we have
\begin{align*}
[a^{n-1}, b^2] [b, a^{-1}]^2 &\equivv [a^{n-1}, b^2] [b^2, a^{-1}] = 
\Psi_{a}^{-1}([a^n,b^2]) \equivv \Psi_{a}^{-1}([a^n,b])^2 \\
&= \left([a^{n-1}, b] [b, a^{-1}] \right)^2 \equivv [a^{n-1}, b]^2 [b, 
a^{-1}]^2.
\end{align*}
Hence the statement holds for all \(n \le -1\). Therefore, the claim is proved.

Secondly, we claim that for all \(n \in \mathbb{Z} \setminus \{0\}\) and \(m 
\in \mathbb{N}\), we have \([a^n, b^m] \equivv [a^n, b]^m\). The proof is by 
induction on~\(m\). For the base case \(m = 1\), there is nothing to prove. 
Next, we prove the inductive step. For \(m \ge 2\), we have
\begin{align*}
[a^n, b]^{-1} [a^n, b^{m+1}] &= \Psi_{b}([a^n,b^m]) \equivv 
\Psi_{b}([a^n,b])^m = \left([a^n, b]^{-1} [a^n, b^2]\right)^m \\
&\equivv [a^n,b]^{-m} [a^n, b^2]^m \equivv [a^n,b]^{-m} [a^n, b]^{2m} \\
&= [a^n, b]^m.
\end{align*}
Therefore, the claim is proved.

Thirdly, we claim that \([a^n, b^m] \equivv [a^n, b]^m\) holds for any 
non-zero \(n,m \in \mathbb{Z}\). It remains to prove the statement for \(m \le 
-1\). The proof is by induction on~\(|m|\). For \(k \geq 1\), we have
\begin{align*}
[a^n,b^{-1}]^{-1} [a^n, b^{k-1}] = \Psi_b^{-1}([a^n,b^k]) \equivv 
\Psi_b^{-1}([a^n, b])^k = [a^n,b^{-1}]^{-k}.
\end{align*}
By setting \(k=2\), we see that the base case (\(m=-1\)) holds. The same line 
proves the inductive step. The claim follows.

Finally, we prove that for all non-zero \(n,m \in \mathbb{Z}\), we have 
\([a^n, b^m] \equivv [a,b]^{nm}\). We set~\(\beta \defeq \sigma_1 \sigma_2 
\sigma_3\). It is easy to see that \(\Psi_\beta(a) = b\) and \(\Psi_\beta(b) = 
a\). We have
\begin{align*}
[a^n, b^m] &= \Psi_\beta([b^n, a^m]) = \Psi_\beta([a^m, b^n])^{-1} \equivv 
\Psi_\beta([a^m, b])^{-n} \\
 &= [b^m, a]^{-n} = [a, b^m]^n \equivv [a, b]^{nm}. \qedhere
\end{align*}
\end{proof}
\medskip

Finally, we are in the position to prove Proposition \ref{Proposition1}. We 
have
\begin{align*}
w_3([a,b]) &= [w_3(a), w_3(b)] = [a, aA_{0,2}a^{-1}] = [a, ab^{-1}a^{-2}] =
a^{-1} a^2ba^{-1} a ab^{-1}a^{-2} \\
&= abab^{-1}a^{-2} = [a^{-1}, b^{-1}] [a^{-2}, b^{-1}]^{-1} \equivv [a,b] 
[a,b]^{-2} = [a,b]^{-1}.
\end{align*}
Since the quotient group \(\Z_3\slash \Bd_3 \cong \pi_3(S^2)\) is infinite 
cyclic, and \([a,b]\) generates it, one has \([a,b]^{-1} \not\equivv [a,b]\). 
Therefore, the image of \(w_3\) under \(r_3\) is a non-identity automorphism. 
The proposition is proved.

\section{Closing remarks}\label{ClosingRemarks}

For the sake of completeness, we calculate the action of the whole generating 
set\footnote{In fact, as \cite[Proposition~4.7]{BNS18} shows, in the case 
\(n=3\), the automorphism \(w_3\) is redundant.}
\begin{align*}
{\rm Aut}_c(\PB_3) \cup \{\Psi_{\sigma_1}, \Psi_{\sigma_2}, \chi_3\} \cup 
\{w_3\}
\end{align*}
of~\({\rm Aut}(\PB_3)\) on the cyclic group~\(\Z_3 \slash \Bd_3 \cong 
\pi_3(\S) \cong \mathbb{Z}\). Firstly, by Lemma~\ref{TrivialityOfAction}, the 
image of \({\rm Aut}_c(\PB_3)\) under \(r_3\) is the trivial subgroup. Next, 
since
\begin{align*}
\chi_3([a,b]) = [a^{-1}, b^{-1}] \equivv [a,b],
\end{align*}
the image of~\(\chi_3\) under~\(r_3\) is the identity automorphism. It remains 
to calculate the action of \(\Psi_{\sigma_1}\) and \(\Psi_{\sigma_2}\). It is 
easy to check that the following identities hold:
\begin{align*}
\sigma_1^{-1} b \sigma_1 &= c, \\
\sigma_2^{-1} a \sigma_2 &= b^{-1} c b.
\end{align*}
We have
\begin{align*}
\Psi_{\sigma_1}([a,b]) &= [\Psi_{\sigma_1}(a), \Psi_{\sigma_1}(b)] = [a, c] = 
a^{-1}c^{-1}ac \overset{(\ref{Conj1})}{=} a^{-1}bab^{-1} = [a,b^{-1}] \equivv 
[a,b]^{-1}
\end{align*}
and
\begin{align*}
\Psi_{\sigma_2}([a,b]) &= [\Psi_{\sigma_2}(a), \Psi_{\sigma_2}(b)] = [b^{-1} c 
b, b] = b^{-1} [c,b] b = b^{-1} c^{-1} b^{-1} c b^2 \\
\overset{(\ref{Conj3})}&{=} b^{-1} (ba b a^{-1} b^{-1})^{-1} b^2 = a b^{-1} 
a^{-1} b = [a^{-1}, b] \equivv [a,b]^{-1}.
\end{align*}
Hence the images of both \(\Psi_{\sigma_1}\) and \(\Psi_{\sigma_2}\) under 
\(r_3\) are non-identity automorphisms.

It seems plausible that for any \(n \ge 3\), the automorphism \(w_n\) acts by 
inversion, and each braid conjugation \(\Psi_\beta\) acts by the sign of the 
permutation \(\pi_\beta\).

\begin{conjecture}\label{Conjecture}
For any \(n \ge 3\), we have \({\rm Im}(r_n) = \{{\rm id}, {\rm -id}\}\). 
Moreover,
\begin{align*}r_n(w_n) = r_n(\Psi_{\sigma_1}) = \ldots = 
r_n(\Psi_{\sigma_{n-1}}) = {\rm -id}.
\end{align*}
\end{conjecture}

This conjecture would imply that the homomorphism \(r_n\) is trivial if and 
only if the exponent of the group \(\pi_n(\S)\) equals two.

\appendix
\section{Proof of Lemma \ref{KeyIsomorphism} and related issues}
\label{AppendixSection}

Firstly, we prove that \(\Bd_n\) is a normal subgroup of \(\Z_n\).
For the inclusion~\(\Bd_n \subseteq \Z_n\),
it is enough to show that~\(\partial_n(\Bd_n) = \{\triv\}\) and \(d_i(\Bd_n) = 
\{\triv\}\) for any \(i \in \{1,2,\ldots,n\}\). These equalities follow from 
Proposition~\ref{BdExpressionProposition}. Next, the subgroup~\(\Bd_n\) is 
normal in~\(\B_n\) by~\cite[Lemma~3.10]{LW09}, and thus,~\(\Bd_n\) is normal 
in~\(\Z_n\). Therefore, the result follows.

At the same time, to prove the inclusion \(\Bd_n \subseteq \Z_n\), the authors 
of \cite{BCWW06} refer to their Lemma~4.1.1 about Delta-groups (see 
Definition~\ref{Delta-groupsDef} below).
To apply this result, they establish Lemma~6.5.2 about the properties of the 
homomorphism~\(\partial_n\). However, part~(3) of Lemma~6.5.2, whose proof is 
omitted, is false. Namely, in our notation,
we have
\[\partial_{n-1} \circ \partial_n \neq \partial_{n-1} \circ d_1\]
for any~\(n \geq 4\) (or, in the notation of~\cite{BCWW06}, \(\partial \circ 
\partial \ne \partial \circ d_0\)). Indeed,
\begin{align}\label{DeltaInequality}
\left(\partial_{n-1} \circ d_1 \right)(A_{1,2}) = \partial_{n-1}(\triv) = 
\triv \neq z_{n-2}^2 \overset{(\ref{PartialA01})}{=} \partial_{n-1}(A_{0,1}) = 
\left(\partial_{n-1} \circ \partial_n\right)(A_{1,2}).
\end{align}
The above inaccuracy causes other issues. 
For instance, part~(1) of Lemma~2.3 in~\cite{LW09} is false. Namely, 
by~(\ref{PartialA01}), in our notation, we have \(\partial_n(A_{0,1}) \neq 
\triv\) (or, in the notation of~\cite{LW09}, \(\partial(A_{0,1}) \neq 
1\)).\footnote{Compare this with~\cite[Lemma~6.5.2]{BCWW06} in which the 
equality~\(\partial(A_{-1,j}) = 1\) is addressed only for~$j\neq 0$.}
As mentioned above, part (3) of~\cite[Lemma~2.3]{LW09}, which is a 
reformulation of part (3) of~\cite[Lemma~6.5.2]{BCWW06}, is false too. Next, 
due to these inaccuracies, part~(2) of~\cite[Lemma~3.2]{LW09} is false: in 
their notation, \(d_0 \chi_{\sigma_1} \neq \theta d_0\) for any~\(n \ge 
3\).\footnote{Note that in~\cite[p.~540]{LW09}, \(d_0\) denotes~\(\partial\).}
Namely, in our notation,
\begin{align*}
\left(\partial_n \circ \Psi_{\sigma_1^{-1}}\right)(A_{1,2}) = 
\partial_n(A_{1,2}) = A_{0,1} \neq \theta_{n-1}(A_{0,1}) = \left(\theta_{n-1} 
\circ \partial_n\right)(A_{1,2})
\end{align*}
for any~\(n \ge 3\), where the inequation holds
because
\[d_1(A_{0,1}) = \triv \neq z_{n-2}^2 = d_1(\theta_{n-1}(A_{0,1}))\]
for \(n \ge 4\) and \(A_{0,1} = A_{1,2}^{-1} \neq A_{1,2} = 
\theta_{2}(A_{0,1})\) for \(n=3\).
The above issues might affect other results of \cite{LW09}. However, we 
verified all the results of \cite{LW09} that we refer to in the present paper 
(namely, Lemma~3.1, Lemma 3.10, Proposition~3.15, and Lemma~4.3).

Secondly, we discuss the isomorphism~\(\Z_n\slash \Bd_n \cong \pi_n(\S)\). To 
prove this result, which is Theorem 1.3 in \cite{BCWW06}, the authors establish 
an isomorphism between a particular Delta-group and Milnor's free group 
construction \(F[S^1]\) on the simplicial one-sphere \(S^1\). Therefore, they 
show that the sequence (\ref{BCWWSequence}) is isomorphic to the Moore complex 
of~\(F[S^1]\). However, their proof contains an inaccuracy, which we fix below.

\begin{figure}[!htb]
\centering
\begin{tikzcd}
{\mathcal G_n} && {\mathcal G_{n-1}} \\
\\
{\mathcal G_{n-1}} && {\mathcal G_{n-2}}
\arrow["{d_i^{\mathcal G}}", from=1-1, to=1-3]
\arrow["{d_j^{\mathcal G}}", from=1-3, to=3-3]
\arrow["{d_{j+1}^{\mathcal G}}"', from=1-1, to=3-1]
\arrow["{d_i^{\mathcal G}}"', from=3-1, to=3-3]
\end{tikzcd}
\caption{Definition of a Delta-group.}\label{DeltaIdentity}
\end{figure}

\begin{definition}\label{Delta-groupsDef}
Let \(\mathcal{G} = (\mathcal{G}_0, \mathcal{G}_1, \ldots)\) be a collection 
of groups. A {\it Delta-group} structure on \(\mathcal{G}\) is a sequence of 
group homomorphisms \(d_0^\mathcal{G}, d_1^\mathcal{G}, \ldots, d_n^\mathcal{G} 
\colon \mathcal{G}_n \to \mathcal{G}_{n-1}\) that satisfy the following 
property: for any~\(n \ge 2\) and any~\(i,j \in \{0,1,\ldots,n-1\}\) such that 
\(i \le j\), one has~\(d_j^\mathcal{G} \circ d_i^\mathcal{G} = 
d_i^\mathcal{G}\circ d_{j+1}^\mathcal{G}\). In other words, the diagram shown 
in Figure \ref{DeltaIdentity} commutes.
\end{definition}
Delta-groups are also referred to as {\it semi-simplicial groups}.

For example, the collection \(\AP \defeq (\PB_1, \PB_2, \PB_2, \ldots)\) 
carries a natural Delta-group structure given by \(d_i^\AP \defeq d_{i+1}\). %
In fact, it can be extended to a so-called simplicial group structure 
(see~\cite[Example~4.7]{CW08}).

Recall that in the present paper, we denote by \(d_i\) (with an empty upper 
index) the removal of the~\(i\)th strand homomorphism, and we number the 
strands from \(1\) to \(n\). In particular, we do not use the plain symbol 
\(d_0\).

Let \(\mathbb{P}_n \defeq \PB_n\), \(d_0^\mathbb{P} \defeq \partial_n\), and 
\(d_i^\mathbb{P} \defeq d_i\) for each~\(i \in \{1,2,\ldots,n\}\).

As stated in~\cite[Corollary~6.5.3]{BCWW06} and~\cite[Proposition~2.4]{LW09}, 
the collection~\(\mathbb{P}\) is a Delta-group. Unfortunately, in this wording, 
this is false. Namely, by~(\ref{DeltaInequality}), we have~\({d_0^\mathbb{P} 
\circ d_0^\mathbb{P} \neq d_0^\mathbb{P} \circ d_1^\mathbb{P}}\). Inspired by 
Proposition~2.9 of~\cite{LW09}, which addresses the fibrantness\footnote{For 
the definition of a fibrant Delta-group, see, e.g., \cite[p.~286]{BCWW06} or 
\cite[p.~534]{LW09}.} of~\(\mathbb{P}\), we suggest the following question on 
the existence of Delta-group structures on pure braid groups.

\begin{question}
Do there exist homomorphisms \(f_n \colon \PB_n \to \PB_{n-1}\) such that the 
collection given by~\(\mathbb{P}_n \defeq \PB_n\), \(d_0^\mathbb{P} \defeq 
f_n\), and \(d_i^\mathbb{P} \defeq d_i\) is a fibrant Delta-group?
\end{question}

To settle the issues mentioned above, we proceed as follows. To prove their 
Theorem~1.3, the authors of~\cite{BCWW06} refer to another collection of groups
\[\mathcal{G}_n \defeq {\rm Ker}(d_{n+1}\colon \PB_{n+1} \to \PB_n)\]
and they set ~\(d_0^\mathcal{G} \defeq \partial_{n+1}\) and \(d_i^\mathcal{G} 
\defeq d_i\) for \(i \in \{1, 2, \ldots, n\}\).\footnote{The collection 
\(\mathcal{G}\) is a version of the Moore loop of \(\mathbb{P}\); cf. 
\cite[p.~533]{LW09} or \cite[p.~288]{BCWW06}.} In fact, to complete the proof 
of Theorem~1.3, it is enough to show that this collection is a Delta-group.

\begin{lemma}
The collection {\normalfont \(\mathcal{G}\)} is a Delta-group.
\end{lemma}
\begin{proof}
We are in the position to prove that
\(d_j^\mathcal{G} \circ d_i^\mathcal{G} = d_i^\mathcal{G}\circ 
d_{j+1}^\mathcal{G}\).

For \(i,j \in \{1,2\ldots,n-1\}\) such that~\(i \le j\),
this follows from the computations in the second line of~\cite[p.~307]{BCWW06}.

For~\(j \in \{1,2,\ldots,n-1\}\) and~\(i = 0\), the result follows from 
part~(2) of~\cite[Lemma~6.5.2]{BCWW06}.

It remains to prove the identity for \((i,j) = (0,0)\), that is, \(\partial_n 
\circ \partial_{n+1} = \partial_n \circ d_1\).
Recall that \(\mathcal{G}_n = \langle A_{k,n+1} \mid 1 \le k \le n\rangle 
\subseteq \PB_{n+1}.\) For \(k=1\), we have
\begin{align*}
\left(\partial_n \circ \partial_{n+1}\right)(A_{1,n+1}) = \partial_n(A_{0, n}) 
= \triv = \partial_n(\triv) = \partial_n(d_1(A_{1,n+1})) = \left(\partial_n 
\circ d_1\right)(A_{1,n+1}),
\end{align*}
and for \(k \in \{2,3,\ldots,n\}\), we have

\[\left(\partial_n \circ \partial_{n+1}\right)(A_{k,n+1}) = \partial_n(A_{k-1, 
n}) = \partial_n(d_1(A_{k,n+1})) = \left(\partial_n \circ 
d_1\right)(A_{k,n+1}). \qedhere\]
\end{proof}
\medskip

Therefore, the proof of Lemma \ref{KeyIsomorphism} is complete.

\end{document}